\documentclass[11pt]{article}
\usepackage{amsfonts}
\usepackage{mathrsfs}
\usepackage{bbm}
\usepackage{cite}
\usepackage{amssymb,amsmath,graphicx}
\usepackage{float} \usepackage[colorlinks=true]{hyperref}
\hypersetup{urlcolor=blue, citecolor=red}

\allowdisplaybreaks

\begin{document}
\title{Boundedness to a logistic chemotaxis system\\
with singular sensitivity}
\author{Xiangdong Zhao{\thanks{
E-mail: zhaoxd1223@163.com (X. D. Zhao) }}\\[6pt]
{\small  School of Mathematics, Liaoning Normal University, Dalian 116029, P.R. China } }
\date{}
\newtheorem{theorem}{Theorem}
\newtheorem{definition}{Definition}[section]
\newtheorem{lemma}{Lemma}[section]
\newtheorem{proposition}{Proposition}[section]
\newtheorem{corollary}{Corollary}[section]
\newtheorem{remark}{Remark}
\renewcommand{\theequation}{\thesection.\arabic{equation}}
\catcode`@=11 \@addtoreset{equation}{section} \catcode`@=12
\maketitle{}
\begin{abstract}

In this paper, we study the parabolic-elliptic Keller-Segel system with singular sensitivity
and logistic-type source:~$
u_t=\Delta u-\chi\nabla\cdot(\frac{u}{v}\nabla v)
+ru-\mu u^k$, $0=\Delta v-v+u$ under the non-flux boundary conditions in a smooth bounded convex domain $\Omega\subset\mathbb{R}^n$,
 $\chi,r,\mu>0$, $k>1$ and $n\ge 2$. It is shown that the system possesses a globally bounded classical solution if $k>\frac{3n-2}{n}$, and $r>\frac{\chi^2}{4}$ for $0<\chi\le 2$, or $r> \chi-1$ for $\chi>2$. In addition, under the same condition for $r,\chi$, the system admits a global generalized solution when $k\in(2-\frac{1}{n},\frac{3n-2}{n}]$, moreover this global generalized solution should be globally bounded provided $\frac{r}{\mu}$ and the initial data $u_0$ suitably small.
\begin{description}
\item[2010MSC:]  35K55; 35B45; 35B40; 92C17
\item[Keywords:] Keller-Segel system;  Singular sensitivity; Logistic source; Boundedness
\end{description}
\end{abstract}

\section{Introduction}
Chemotaxis, is a spontaneous cross-diffusion phenomena by which organisms direct their movements in regard to a
stimulating chemical. In 1970, Keller and Segel proposed a model to represent the chemotaxis phenomena, i.e., the oriented or partially oriented movement of cells with respect to a chemical signal produced by the cells themselves \cite{EL}:
\begin{eqnarray}\label{p1}
 \left\{
\begin{array}{llll}
u_t=\Delta u-\chi\nabla\cdot(\frac{u}{v}\nabla v)
, & x\in\Omega,~t>0,\\[4pt]
\displaystyle \tau  v_t=\Delta v-v+u,& x\in\Omega,~t>0,
\end{array}\right. \end{eqnarray}
where $\tau\in\{0,1\}$,$\chi>0$. The singular chemotactic sensitive function $\frac{\chi}{v}$ with $\chi>0$ is derived by the Weber-Fechner law on the response of the cells $u$ to the stimulating chemical signal $v$.  With the singularity determined by the sensitive function $\frac{\chi}{v}$, the cellular movements are governed by the taxis flux $\frac{\chi\nabla v}{v}$, which may be unbounded when $v \approx 0$. Different to the classical Keller-Segel model (i.e., replacing the singular sensitive function $\frac{\chi}{v}$ by the constant function $\chi$ in \eqref{p1}), it is important to obtain a lower bound on $v$ for studying the global dynamical behavior. This can be achieved by a pointwise estimate \cite{KMT}
\begin{align*}
v(x,t)\ge c_0 \int_\Omega u(x,t)dx, ~~~~~~~x\in\Omega,~t>0
\end{align*}
with $c_0=c_0(|\Omega|,n)>0$. Due to the mass conservation of cells $u$ in system \eqref{p1}, it is known that the singularity involved in sensitive function $\frac{\chi}{v}$ is in fact absent. Generally, chemotactic sensitive coefficient $\chi>0$ properly small benefits the global existence-boundedness of solutions to system \eqref{p1}, which can be presented in \cite{KOA,K,B,LW,FST,FS}. It is pointed that for the parabolic-elliptic case of the system \eqref{p1} $(\tau=0)$ with radial assumption, Nagai and Senba proved that the problem admits a finite time blow-up solution \cite{TT} if $\chi>\frac{2n}{n-2}$ with $n\ge 3$, and $\int_\Omega u_0|x|^2dx$ sufficiently small.

 Consider the chemotaxis system as follows
\begin{eqnarray}\label{p2}
 \left\{
\begin{array}{llll}
u_t=\Delta u-\chi\nabla\cdot(u\nabla v)+ru-\mu u^k, & x\in\Omega,~t>0,\\[4pt]
\displaystyle  \tau v_t=\Delta v-v+u,& x\in\Omega,~t>0,\\[4pt]
\end{array}\right. \end{eqnarray}
where $\chi,r,\mu>0$, $k>1$ and $\tau\in\{0,1\}$. Such self-limiting growth mechanism involved in the logistic-type source generally benefits the global dynamic of solutions. For parabolic-elliptic case of \eqref{p2} ($\tau=0$ ), the system with $k=2$ possesses a global weak solution if $\mu>0$ and a global bounded classical solution if $\mu>\frac{n-2}{n}\chi$ \cite{JM}.
If $k>2-\frac{1}{n}$ with $n\ge 1$, there exists a global very weak solution, which is globally bounded provided  $\mu$ sufficiently large and $u_0$ sufficiently small  \cite{M3}.  Replacing $0=\Delta v-v+u$ in \eqref{p2} by $0=\Delta v-m(t)+u$ with $m(t):=\frac{1}{|\Omega|}\int_\Omega u(x,t)dx$, it is shown with radial assumption that the system admits a finite time blow-up solution if $1<k<\frac{3}{2}+\frac{1}{2n-2} $ with $n\ge 5$ \cite{M2}. For the parabolic-parabolic case of \eqref{p2} ($\tau=1$), if $k=2$, $n=2$ \cite{OO}, or $n\ge 3$ with $\mu>0$ sufficiently large \cite{W1}, the problem possesses globally bounded classical solutions. If $k>2-\frac{1}{n}$ with $n\ge 1$, there exists global very weak solutions \cite{V}, which are globally bounded provided $\frac{r}{\mu}$ and the initial data all sufficiently small for $n=3$ \cite{V1}. In addition, more properties of solutions to \eqref{p2} can be found in \cite{L1,TM,M4}.

Recall the chemotaxis system with singular sensitivity and logistic source
 \begin{eqnarray}\label{p22}
 \left\{
\begin{array}{llll}
u_t=\Delta u-\chi\nabla\cdot(\frac{u}{v}\nabla v)+ru-\mu u^k, & x\in\Omega,~t>0,\\[4pt]
\displaystyle  \tau v_t=\Delta v-v+u,& x\in\Omega,~t>0,
\end{array}\right. \end{eqnarray}
where $\chi,r,\mu>0$, $k>1$ and $\tau\in\{0,1\}$. The difficulty in studying global solvability of solutions comes from the hazardous combination of singular sensitive chemotactic mechanism and the self-limiting growth mechanism involved in logistic source. Due to missing the mass conservation for $u$, the singularity contained in chemotactic term may be presence. Similarity to system \eqref{p1}, the suitable smallness of chemotactic sensitive coefficient $\chi>0$ is necessary to establish global existence-boundedness of solutions to system \eqref{p22}.  If $n,k=2$, there exists a unique globally bounded classical solution \cite{FF,ZZ}, whenever
 \begin{equation}\label{th21}
r>
\begin{cases}
    &\frac{\chi^2}{4},~~~~  0<\chi\le 2,\\
    &\chi-1 ,~~~ ~\chi>2.
    \end{cases}
\end{equation}
In addition, the author has proved for $k>2-\frac{1}{n}$ that the system \eqref{p22} with $\tau=1$ possesses a global very weak solution provided  $\chi$ suitably small related to $r,k$ \cite{ZZZ}.

Turn to a chemotaxis-consumption model of the type
\begin{eqnarray}\label{p3}
 \left\{
\begin{array}{llll}
u_t=\Delta u-\chi\nabla\cdot(\frac{u}{v}\nabla v)+ru-\mu u^k, & x\in\Omega,~t>0,\\[4pt]
\displaystyle  v_t=\Delta v-uv,& x\in\Omega,~t>0,\\[4pt]
\end{array}\right.
\end{eqnarray}
where $\chi,r,\mu>0$ and $k>1$.
It seems difficult to study the global existence and further global dynamic behavior of solutions to \eqref{p3} than that for the problems \eqref{p2} and \eqref{p22}.  The troubles lie in the interplay of the consumptive effect
with the singular chemotactic mechanism and self-limiting logistic source.  Intuitively, the oxygen $v$ shrinks in \eqref{p3}$_2$
during evolution, and then enhances the chemotactic strength of the bacteria because of the
singular behavior when $v\approx 0$ in \eqref{p3}$_1$. This implies the singularity in chemotactic sensitive function $\frac{\chi}{v}$ should be persistence.
Recall some results on the case of  $k=2$. It is shown for $n\ge 2$ that there exists a global classical solution if $0<\chi<\sqrt{\frac{2}{n}}$ and $\mu>\frac{n-2}{2n}$, and that for $n=1$ the global classical solution is globally bounded if $\chi,r,\mu>0$ \cite{LL1}. Moreover, the problem for $n\ge 1$ possesses a global generalized solution \cite{LL}. If $k>1+\frac{n}{2}$, the author has established the global solvability of classical solutions \cite{Z1}.
 We refer to \cite{M5,M6,M7} for more results on chemotaxis-consumption system without logistic source.

 In this paper, we consider the following parabolic-elliptic chemotaxis system with singular sensitivity and logistic-type source
\begin{eqnarray}\label{p}
 \left\{
\begin{array}{llll}
u_t=\Delta u-\chi\nabla\cdot(\frac{u}{v}\nabla v)
+ru-\mu u^k, & x\in\Omega,~t>0,\\[4pt]
\displaystyle  0=\Delta v-v+u,& x\in\Omega,~t>0,\\[4pt]
  \displaystyle \frac{\partial u}{\partial {\nu}}=\frac{\partial v}
  {\partial {\nu}}=0 ,& x\in\partial\Omega,~t>0,\\[4pt]
  \displaystyle u(x,0)=u_0(x),  &x\in\Omega,
\end{array}\right. \end{eqnarray}
where $\chi,r,\mu>0$, $k>1$. $\Omega\subset\mathbb{R}^n$ $(n\ge 2)$ is a smooth bounded convex domain, $\frac{\partial}{\partial \nu}$ denotes the derivation with respect to the outer normal of $\partial\Omega$, and the initial data
\begin{align}\label{con}
u_0(x)\in C^0(\overline{\Omega}),~ u_0(x)\ge 0 {\rm ~and~} u_0(x)\neq 0,~ x\in \overline\Omega.
\end{align}

To study the global dynamic behavior of solution to system \eqref{p} for the more general exponent $k>1$ in the logistic-type source $ru-\mu u^k$, we introduce the generalized solution to \eqref{p} via the following definitions inspired by \cite{M3,LL}.

\begin{definition}\label{def1} {\rm
 A pair $(u,v)$ of nonnegative functions
$$u\in L_{\rm loc}^1(\Omega\times(0,\infty)),~v\in L_{\rm loc}^1((0,\infty);W^{1,1}(\Omega))$$
will be called a {\em very weak subsolution} to \eqref{p}, if
\begin{align*}
u^k~{\rm and} ~\frac{u}{v}\nabla v ~{\rm belong ~to}~L_{\rm loc}^1(\Omega\times(0,\infty)),
\end{align*}
and moreover
\begin{align}\label{def11}
&-\int_0^\infty\int_\Omega u\varphi_t -\int_\Omega u_0\varphi(\cdot,0)\le \int_0^\infty\int_\Omega u\Delta \varphi+\chi\int_0^\infty\int_\Omega \frac{u}{v}\nabla v\cdot\nabla\varphi\nonumber\\
&~~~~~~~~~~~~~~~~~~~~~~~~~~~~~~~~~~~~~~~~+r\int_0^\infty\int_\Omega u\varphi-\mu\int_0^\infty\int_\Omega u^k\varphi,\\\label{def12}
&-\int_0^\infty\int_\Omega v\psi_t-\int_\Omega v_0\psi(\cdot,0)+\int_0^\infty\int_\Omega \nabla v\cdot\nabla \psi+\int_0^\infty\int_\Omega v\psi=\int_0^\infty\int_\Omega u\psi
\end{align}
hold for all
\begin{align}\label{def13}
\varphi  &\in C_0^\infty(\bar\Omega\times(0,\infty))~{\rm with}~\varphi\ge 0~{\rm and}~\frac{\partial \varphi}{\partial \nu}=0~{\rm on}~\partial\Omega\times(0,\infty),\\\label{def14}
\psi &\in L^\infty(\Omega\times(0,\infty))\cap L^2((0,\infty);W^{1,2}(\Omega)), ~{\rm and}~\psi_t\in L^2(\Omega\times(0,\infty)).
\end{align}}
\end{definition}

\begin{definition}\label{def2}{\rm
A pair of nonnegative functions $u\in L_{\rm loc}^1(\Omega\times(0,\infty)),~v\in L_{\rm loc}^1((0,\infty);W^{1,1}(\Omega))$
form a {\em weak logarithmic supersolution} to \eqref{p}, if
$$\frac{u^k}{1+u}, \frac{|\nabla u|^2}{(1+u)^2}~{\rm and}~ \frac{|\nabla v|^2}{ v^{2}}~{\rm belong ~to}~L_{\rm loc}^1(\Omega\times(0,\infty)),$$
and moreover
\begin{align}\label{def21}
-\int_0^\infty\int_\Omega \ln(1+u)\varphi_t-\int_\Omega \ln(1+u_0)\varphi(\cdot,0)&\ge\int_0^\infty\int_\Omega \frac{|\nabla u|^2}{(1+u)^2}\varphi-\chi\int_0^\infty\int_\Omega \frac{u}{(1+u)^2v}\nabla u\cdot\nabla v \varphi\nonumber\\
&~~-\int_0^\infty\int_\Omega\frac{\nabla u \cdot\nabla \varphi}{1+ u}+\chi\int_0^\infty\int_\Omega\frac{u\nabla v}{(1+u)v} \cdot\nabla\varphi\nonumber\\
&~~
+r\int_0^\infty\int_\Omega \frac{u}{1+u}\varphi-\mu\int_0^\infty\int_\Omega \frac{u^k}{1+u}\varphi
\end{align}
and the equality  \eqref{def12} hold for all $\varphi$ and $\psi$ satisfying \eqref{def13} and \eqref{def14}.}
\end{definition}

\begin{definition}\label{def3}{\rm
A couple of function $(u,v)$ will be called a {\em generalized solution} to \eqref{p} if it is both a very weak subsolution and a weak logarithmic supersolution of \eqref{p}. } \end{definition}

To obtain the global dynamic behavior of solution to system \eqref{p} for general $k>1$, we will at first establish a positive uniform-in-time lower bound for chemical signal $v$. With the aid of a crucial ODE inequality \cite{SSW}, this will be realized by a uniform-in-time upper estimate on the integral $\int_\Omega u^{-m}dx$ with some $m>0$. Furthermore, via a standard process on energy estimate for $\int_\Omega u^pdx$ with $p>1$,  we conclude the global boundedness of the classical solutions if $k>\frac{3n-2}{n}$ and $\chi>0$ suitably small relative to $r>0$.

 In order to deal with the global generalized solution of classical parabolic-elliptic system \eqref{p2}(i.e., without singular sensitive function $\frac{1}{v}$ in \eqref{p}) \cite{V}, the crucial step is to conclude the relative compactness of the solution $\{v_\epsilon\}_{\epsilon\in(0,1)}$ to the corresponding  regularization problem in $L_{\rm loc}^\frac{k}{k-1}((0,\infty);W^{1,\frac{k}{k-1}}(\Omega))$ with respect to the strong topology for $k\in(2-\frac{1}{n},2)$. Differently, for the system \eqref{p}, we will firstly show that $\{v_\epsilon\}_{\epsilon\in(0,1)}$ has a uniform-in-time lower bound (indepedent of $\epsilon\in(0,1)$).
Secondly, it will be derived that for some $p>\frac{nk}{n-1}$ with $k>2-\frac{1}{n}$ that
$\{\frac{\nabla v_\epsilon}{v_\epsilon}\Big\}_{\epsilon\in(0,1)}$  is relatively compact in $L_{\rm loc}^p(\Omega\times(0,\infty))$ with respect to the weak topology.
Finally, upon selecting a suitable subsequence, we will obtain a global generalized solution for $k>2-\frac{1}{n}$ with $n\ge 2$ by a standard compactness argument. Furthermore, if $\frac{r}{\mu}$ and the initial data $u_0$ suitably small, this global generalized solution is in fact globally bounded.

Now, we state the main results of this paper.
\begin{theorem}\label{th1}
Let $n\ge 2$ and $k>\frac{3n-2}{n}$. If $r,\chi>0$ satisfying
\begin{equation}\label{th11}
r>
\begin{cases}
    \frac{\chi^2}{4},& 0<\chi\le 2,\\
    \chi-1,&\chi>2,
    \end{cases}
\end{equation}
the problem \eqref{p} possesses a unique globally bounded classical solution.
\end{theorem}

\begin{theorem}\label{th2}
Let $n\ge 2$, $k>2-\frac{1}{n}$ and $r,\chi>0$ satisfy \eqref{th11}. Then the system \eqref{p} admits at least a global generalized solution.
\end{theorem}

\begin{theorem}\label{th3}
Let $(u,v)$ be the global generalized solution for $k\in(2-\frac{1}{n},\frac{3n-2}{n}]$ established in Theorem \ref{th2}. Then for $p>\frac{n(n+2)}{2(n+1)}$ there exist $\eta,\lambda>0$ small such that this solution is globally bounded provided $\frac{r}{\mu}<\eta$ and $\int_\Omega u_0^pdx<\lambda$.
\end{theorem}

\begin{remark}{\rm Since $\frac{3n-2}{n}=2$ for $n=2$, we conclude by Theorem \ref{th1} with \cite{FF} that the classical solution to \eqref{p} for the case $n=2$ must be globally bounded if $k\ge 2$. In addition, Theorems \ref{th2} and \ref{th3} show that $k<2$  is permitted for the global existence-boundedness of solution to \eqref{p}. This extends the global boundedness results for \eqref{p} with $k=2$ obtained in \cite{FF}.} \end{remark}

\begin{remark}\label{re1}
{\rm Recall from \cite{JM,GST,M3} that the classical parabolic-elliptic chemotaxis system \eqref{p2} possesses a global or further globally bounded classical (or generalized) solution if $k>1$ and $\mu>0$ suitably large in logistic source $ru-\mu u^k$. While Theorems above say that the global solvability of solution to the system \eqref{p} requires not only the restriction on logistic kinetics but also the chemotactic sensitive coefficient $\chi>0$ properly small relative to $r>0$, and the same were observed for the parabolic-parabolic case of system 
\eqref{p22} \cite{ZZZ}. Here the difficulty due to the singular sensitivity is substantial. }
\end{remark}

The rest part of the paper is arranged as follows. In Section 2, we will establish a uniform-in-time lower bound estimate on $v$ and demonstrate the global boundedness of the classical solution. Section 3 deals with the global existence of classical solution to the corresponding regularization problems. Then we prove the global existence and boundedness to the generalized  solution to system \eqref{p} in Section 4.

\section{Global boundedness of classical solutions}
We at first give a lemma on the local existence of classical solutions to system \eqref{p} without proof, which can be obtained by the contraction argument as that in \cite[ Lemma 2.1]{FF}.
\begin{lemma}\label{lemma1}
Assume that $u_0$ satisfies \eqref{con}. If $k>1$, $r,\chi,\mu>0$, then
there exist $T_{\max}\in(0,+\infty]$ and a unique pair $(u,v)$ of functions
 \begin{equation*}
   \left\{
    \begin{aligned}
    &u\in C^0(\overline\Omega\times[0,T_{\max}))\cap C^{2,1}
     (\overline{\Omega}\times(0,T_{\max})),\\[4pt]
    &v\in  C^{2,0}
     (\overline{\Omega}\times(0,T_{\max})),
    \end{aligned}
    \right.
\end{equation*}
fulfilling \eqref{p} in the classical sense with $u,v>0$ in $\overline\Omega\times(0,T_{\max})$. Moreover, either $ T_{\max}=\infty$, or ${\limsup}_{t\rightarrow T_{\max}}\| u(\cdot,t)\|_{L^\infty(\Omega)}=\infty$, or ${\liminf}_{t\rightarrow T_{\max}} \inf_{x\in\Omega} v(x,t)=0$.\qquad$\Box$
\end{lemma}

Let $(u,v)$ be the local classical solution in this section. Without loss of generality, assume that $T_{\max}>1$. Then we have the following a priori estimates.
\begin{lemma}\label{lemma2}
For $k>1$ and $r,\chi,\mu>0$, it holds that
\begin{align}
&\int_\Omega udx\le m^*,~
t\in(0,T_{\max}),\label{lem21}\\
&\int_t^{t+1}\int_\Omega u^kdxds\le M_1~~{\rm for~all}~t\in(0,T_{\max}-1)\label{lem23}
\end{align}
with $m^*=\max\Big\{\int_\Omega u_0dx,|\Omega|\Big(\frac{r}{\mu}\Big)^{\frac{1}{k-1}}\Big\}$ and $M_1=\frac{(1+r) m^*}{\mu}$.
\end{lemma}
{\bf Proof.}\
Integrate \eqref{p}$_1$ over $\Omega$ to know
\begin{align}\label{21}
\frac{d}{dt}\int_\Omega udx&=r\int_\Omega udx-\mu\int_\Omega u^kdx\\
&\le r\int_\Omega udx-\frac{\mu}{|\Omega|^{k-1}}\Big(\int_\Omega udx\Big)^{k}\label{22},~t\in(0,T_{\max})
\end{align}
by the H\"{o}lder inequality. We get \eqref{lem21} by the Bernoulli inequality with \eqref{22}.
The estimate \eqref{lem23} comes directly from by integrating \eqref{21} from $t$ to $t+1$.
\qquad$\Box$\medskip

To study the dynamic behavior of solutions to \eqref{p} for $k>1$, we should pay attention to establish a positive uniform-in-time lower bound for $v$. With the aid of the following crucial ODE inequality \cite{SSW}, this will be realized by a uniform-in-time upper estimate on the integral $\int_\Omega u^{-m}dx$ with some $m>0$ \cite{FF}.
\begin{lemma}\cite[Lemma 3.4]{SSW}\label{lemma4}
Let $T>0$, and suppose that $y$ is a nonnegative absolutely continuous function on $[0,T)$ satisfying
\begin{align*}
y'(t)+ay(t)\le f(t)~~~{\rm for~almost~every~}t\in(0,T)
\end{align*}
with some $a>0$ and a nonnegative function $f\in L_{\rm loc}^1([0,T))$ for which there exists $b>0$ such that
\begin{align*}
\int_t^{t+1}f(s)dx\le b~~~{\rm for~all}~t\in[0,T-1).
\end{align*}
Then
\begin{align*}
y(t)\le \max\{y(0)+b,\frac{b}{a}+2b\}~~~{\rm for ~all}~t\in(0,T).
\end{align*}
\end{lemma}

Now, we have:
\begin{lemma}\label{lemma5}
Let $k>1$, $\mu>0$ and $r,\chi>0$ satisfy
 \begin{equation}\label{lem51}
r>
\begin{cases}
    &\frac{\chi^2}{4},~~~~  0<\chi\le 2,\\
    &\chi-1 ,~~~ ~\chi>2.
    \end{cases}
\end{equation}
Then there exists some $\delta_0>0$ such that
\begin{align}\label{lem52}
v(x,t)\ge \delta_0~~~{\rm for~all}~(x,t)\in\Omega\times(0,T_{\max}).
\end{align}
\end{lemma}
{\bf Proof.}\
Since $u\in C^0(\bar{\Omega}\times[0,T_{\max}))$ by Lemma \ref{lemma1}, we know that there exists $t_0\in(0,T_{\max})$ such that
 \begin{align}\label{30}
 \int_\Omega u(x,t)dx\ge \frac{1}{2}\int_\Omega u_0dx,~~t\in(0,t_0].
 \end{align}
Invoke the pointwise lower bound estimate in \cite{KMT} to know
\begin{align}\label{310}
v(x,t)&\ge c_1\int_\Omega u(x,t)dx ~~{\rm ~ for ~all}~ x\in\Omega~{\rm and}~t\in(0,T_{\max})
\end{align}
with some $c_1>0$. Hence, \eqref{310} with \eqref{30} yields
\begin{align}\label{311}
v(x,t)&\ge c_1\int_\Omega u(x,t)dx\ge \frac{c_1}{2}\int_\Omega u_0dx:=\beta_0 ~~{\rm ~ for ~all}~ x\in\Omega~{\rm and}~t\in(0,t_0].
\end{align}

For $m>0$, it is known from  \eqref{p}$_1$ that
\begin{align}\label{31}
\frac{1}{m}\frac{d}{dt}\int_\Omega u^{-m}dx&=-\int_\Omega u^{-m-1}[\Delta u-\chi\nabla \cdot(\frac{u}{v}\nabla v)+ru-\mu u^k]dx\nonumber\\
&=-(m+1)\int_\Omega u^{-m-2}|\nabla u|^2dx+\chi(m+1)\int_\Omega u^{-m-1}\frac{\nabla u\cdot\nabla v}{v}dx\nonumber\\
&~~-r\int_\Omega u^{-m}dx+\mu \int_\Omega u^{-m-1+k}dx,~t\in(t_0,T_{\max}).
\end{align}
If $a\in(0,\chi)$, we get from \eqref{p}$_2$ that
\begin{align}\label{32}
(m+1)a\int_\Omega u^{-m-1}\frac{\nabla u\cdot\nabla v}{v}dx&=-\frac{(m+1)a}{m}\int_\Omega \nabla u^{-m}\cdot\frac{\nabla v}{v}dx\nonumber\\
&\le-\frac{(m+1)a}{m}\int_\Omega u^{-m}\frac{|\nabla v|^2}{v^2}dx+\frac{(m+1)a}{m}\int_\Omega u^{-m}dx,
\end{align}
and by Young's inequality that
\begin{align}\label{33}
(m+1)(\chi-a)\int_\Omega u^{-m-1}\frac{\nabla u\cdot\nabla v}{v}dx&\le (m+1)\int_\Omega u^{-m-2}|\nabla u|^2dx\nonumber\\
&~~+\frac{(m+1)(\chi-a)^2}{4}\int_\Omega u^{-m}\frac{|\nabla v|^2}{v^2}dx.
\end{align}
Let $m:=\frac{4a}{(\chi-a)^2}$ for $a\in(0,\chi)$. Then $\frac{(m+1)(\chi-a)^2}{4}=\frac{(m+1)a}{m}$. Combining \eqref{31}--\eqref{33}, we have
\begin{align}\label{34}
\frac{1}{m}\frac{d}{dt}\int_\Omega u^{-m}dx&\le -(r-\frac{(m+1)a}{m})\int_\Omega u^{-m}dx+\mu\int_\Omega u^{-m-1+k}dx,~t\in(t_0,T_{\max}).
\end{align}
Denote
$$f(a):=-4(r-\frac{(m+1)a}{m})=a^2-(2\chi-4)a+\chi^2-4r.$$
A direct calculation shows that  $\Delta=16(r+1-\chi)>0$ for $r>\max\{\chi-1,0\}$, and hence $f(a)<0$ for $a\in(a_{-},a_{+})$,
here $a_{\pm}=\chi-2\pm2\sqrt{r+1-\chi}$.
By the Vi\`{e}te formula, we know  $a_-<0<a_+$ if $r>\frac{\chi^2}{4}$ with $\chi>0$, and $0<a_{-}<\chi<a_{+}$ if $\chi-1<r\le \frac{\chi^2}{4}$ with $\chi>2$. Therefore, if $r,\chi>0$ satisfying \eqref{lem51}, there exists some $c_0>0$ such that $-(r-\frac{(m+1)a}{m})=\frac{f(p)}{4}\le -c_0<0$ for $a\in(0,\chi)\cap(a_{-},a_{+})$.

If $k-1-m=0$, it is known from \eqref{34} that
\begin{align}\label{35}
\frac{1}{m}\frac{d}{dt}\int_\Omega u^{-m}dx&\le -c_0\int_\Omega u^{-m}dx+\mu|\Omega|,~t\in(t_0,T_{\max}).
\end{align}

If $k-1-m<0$, we obtain by Young's inequality with \eqref{34} that
\begin{align}\label{36}
\frac{1}{m}\frac{d}{dt}\int_\Omega u^{-m}dx&\le -\frac{c_0}{2}\int_\Omega u^{-m}dx+\mu|\Omega|(\frac{2\mu}{c_0})^\frac{m+1-k}{k-1},~t\in(t_0,T_{\max}).
\end{align}

Similar process for the case of $k-1-m>0$, we get
\begin{align}\label{37}
\frac{1}{m}\frac{d}{dt}\int_\Omega u^{-m}dx&\le -c_0\int_\Omega u^{-m}dx+\int_\Omega u^kdx+\mu^\frac{k}{m+1}|\Omega|,~t\in(t_0,T_{\max}).
\end{align}
The estimates \eqref{35}--\eqref{37} show for $k>1$ that
\begin{align}\label{38}
\frac{1}{m}\frac{d}{dt}\int_\Omega u^{-m}dx&\le -\frac{c_0}{2}\int_\Omega u^{-m}dx+\int_\Omega u^kdx+C_1,~t\in(t_0,T_{\max}).
\end{align}
with $C_1=\mu|\Omega|(\frac{2\mu}{c_0})^\frac{m+1-k}{k-1}+\mu^\frac{k}{m+1}|\Omega|$. Based on Lemma \ref{lemma4} with \eqref{lem23} and \eqref{38}, then
\begin{align}\label{39}
\int_\Omega u^{-m}dx\le C_2, ~t\in(t_0,T_{\max})
\end{align}
with $C_2=\{ (C_1+M_1)m+m\int_\Omega u(x,t_0)^{-m}dx,~(C_1+M_1)m+\frac{2m^2(C_1+M_1)}{c_0}\}$.

Let $\alpha:=\frac{m}{m+1}\in (0,1)$. Then for $k>1$ we obtain from \eqref{39} and \eqref{310} that
\begin{align*}
v(x,t)\ge c_1\int_\Omega udx\ge c_1|\Omega|^\frac{m+1}{m}\Big(\int_\Omega u^{-m}dx\Big)^{-\frac{1}{m}}\ge c_1C_2^{-\frac{1}{m}}|\Omega|^\frac{m+1}{m} =:\eta_0
\end{align*}
by the H\"{o}lder inequality for all $(x,t)\in\Omega\times(t_0,T_{\max})$. This together with \eqref{311} concludes \eqref{lem52} with $\delta_0=\min\{\beta_0,\eta_0\}$.  \qquad$\Box$\medskip

Based on the uniform-in-time lower bound estimate for $v$ in Lemma \ref{lemma5}, we establish the following $L^p$-estimate for $u$.
\begin{lemma}\label{lemma25}
If $k>\frac{3n-2}{n}$, $\mu>0$ and $r,\chi>0$ satisfying \eqref{lem51}, then for $p>1$ there exists some $M_2>0$ such that
\begin{align}\label{lem251}
\int_\Omega u^pdx\le M_2,~t\in(0,T_{\max}).
\end{align}
\end{lemma}
{\bf Proof.}\
A simple calculation with \eqref{p}$_1$ and \eqref{lem52} shows
\begin{align}\label{251}
\frac{1}{p}\frac{d}{dt}\int_\Omega u^pdx&=-(p-1)\int_\Omega u^{p-2}|\nabla u|^2dx+\chi(p-1)\int_\Omega \frac{u^{p-1}}{v}\nabla u\cdot\nabla vdx\nonumber\\
&~~~+r\int_\Omega u^pdx-\mu\int_\Omega u^{p+k-1}dx\nonumber\\
&\le -\frac{1}{p}\int_\Omega u^pdx+\frac{\chi^2(p-1)}{4\delta_0^2}\int_\Omega u^{p}|\nabla v|^2dx-\frac{\mu}{2}\int_\Omega u^{p+k-1}dx+C_3\nonumber\\
&\le -\frac{1}{p}\int_\Omega u^pdx+C_4\int_\Omega |\nabla v|^{\frac{2(p+k-1)}{k-1}}dx-\frac{\mu}{4}\int_\Omega u^{p+k-1}dx+C_3
\end{align}
with $C_3=(r+1)(\frac{2(r+1)}{\mu})^\frac{p}{k-1}$ and $C_4=\frac{\chi^2(p-1)}{4\delta_0^2}(\frac{\chi^2(p-1)}{\delta_0^2\mu})^\frac{2p}{k-1}$ for $t\in(0,T_{\max})$. Since $\int_\Omega v(x,t)dx=\int_\Omega u(x,t)dx\le m^*$, $t\in(0,T_{\max})$, invoking the classical result by Br\'{e}zis and Strauss \cite{BS} and the Minkowski inequality, we get for $r\in(1,\frac{n}{n-1})$ that
\begin{align}\label{252}
\|v\|_{W^{1,r}(\Omega)}\le C_{BS}\|\Delta v\|_{L^1(\Omega)}\le C_{BS}\|v-u\|_{L^1(\Omega)}\le 2C_{BS}m^*
\end{align}
 with some $C_{BS}>0$. According to the standard elliptic $L^p$-theory,  we know from \eqref{p}$_2$ for $m\ge 1$ that
 \begin{align}\label{2521}
 \|v\|_{W^{2,m}(\Omega)}\le C_5\|u\|_{L^{m}(\Omega)}
 \end{align}
  with some $C_5>0$, and thus by the Gaglirado-Nirenberg inequality with \eqref{252} that
\begin{align}\label{253}
\|\nabla v\|_{L^{\frac{2(p+k-1)}{k-1}}(\Omega)}\le C_{GN}\| v\|_{W^{2,p+k-1}(\Omega)}^a\|\nabla v\|_{L^r(\Omega)}^{1-a}\le C_6\| u\|_{L^{p+k-1}(\Omega)}^a
\end{align}
with some $C_6>0$,
 where $$a=\frac{\frac{n}{r}-\frac{(k-1)n}{2(p+k-1)}}{1-\frac{n}{p+k-1}+\frac{n}{r}}.$$
If $k>\frac{3n-2}{n}$, then $a\in(0,1)$ and $\frac{2a}{k-1}<1$. Consequently, we obtain from \eqref{251} and \eqref{253} that
\begin{align*}
\frac{1}{p}\frac{d}{dt}\int_\Omega u^pdx\le -\frac{1}{p}\int_\Omega u^pdx+C_7,~t\in(0,T_{\max}),
\end{align*}
by Young's inequality with $C_7=C_3+C_4|\Omega|(\frac{4C_4}{\mu})^\frac{2a}{k-1-2a}C_6^\frac{2(p+k-1)}{k-1-2a}$. This concludes \eqref{lem251} by the Bernoulli inequality with some $M_2>0$.
 \qquad$\Box$\medskip

{\bf Proof of the Theorem \ref{th1}}\
By the variation-of-constants formula for $u$ and the order preserving of the Neumann heat semigroup $\{{\rm e}^{t\Delta}\}_{t\ge 0}$ with the positivity of $u$, we know
\begin{align}\label{tt1}
u(x,t)&={\rm e}^{t\Delta}u_0-\chi\int_0^t{\rm e}^{(t-s)\Delta}\nabla\cdot(\frac{u}{v}\nabla v)ds+\int_0^t{\rm e}^{(t-s)\Delta}(ru-\mu u^k)ds, \nonumber\\
&\le {\rm e}^{t\Delta}u_0-\chi\int_0^t{\rm e}^{(t-s)\Delta}\nabla\cdot(\frac{u}{v}\nabla v)ds+r\int_0^t{\rm e}^{(t-s)\Delta}uds,~(x,t)\in\Omega\times (0,T_{\max}).
\end{align}
Let $p>n$ in Lemma \ref{lemma25}. Then we have by the classical elliptic eqution theory to \eqref{p}$_2$ with \eqref{2521} and \eqref{lem251} that
\begin{align}\label{tt2}
\|\nabla v\|_{L^\infty(\Omega)}\le C_9\|v\|_{W^{2,p}(\Omega)} \le C_9C_5\|u\|_{L^p(\Omega)}\le C_{10}, ~t\in(0,T_{\max})
\end{align}
with some $C_9,C_{10}>0$. Consequently, invoking the homogeneous Neumann semigroup estimates in \cite[Lemma 1.3]{M1}  with \eqref{lem21}, it is known from \eqref{tt1} with \eqref{lem52}, \eqref{lem251} and  \eqref{tt2} that
\begin{align}\label{tt3}
\|u\|_{L^{\infty}(\Omega)}&\le\|{\rm e}^{t\Delta}u_0\|_{L^{\infty}(\Omega)}+\chi\int_0^t\|{\rm e}^{(t-s)\Delta}\nabla\cdot(\frac{u}{v}\nabla v)\|_{L^{\infty}(\Omega)}ds+r\int_0^t\|{\rm e}^{(t-s)\Delta}u\|_{L^{\infty}(\Omega)}ds\nonumber\\
&\le \|u_0\|_{L^{\infty}(\Omega)}+\frac{\chi K_4}{\delta_0}\int_0^t(1+(t-s)^{-\frac{1}{2}-\frac{1}{p}}){\rm e}^{-\lambda_1(t-s)}\|u\nabla v\|_{L^{p}(\Omega)}ds\nonumber\\
&~~+rK_1\int_0^t(1+(t-s)^{-\frac{1}{2}}){\rm e}^{-\lambda_1(t-s)}\|u-\overline{u}\|_{L^{2}(\Omega)}ds+\frac{rm^*}{|\Omega|}\nonumber\\
&\le \frac{\chi K_4}{\delta_0}\int_0^t(1+(t-s)^{-\frac{1}{2}-\frac{1}{p}}){\rm e}^{-\lambda_1(t-s)}\|u\|_{L^{p}(\Omega)}\|\nabla v\|_{L^\infty(\Omega)}ds\nonumber\\
&~~+2rK_1\int_0^t(1+(t-s)^{-\frac{1}{2}}){\rm e}^{-\lambda_1(t-s)}\|u\|_{L^{2}(\Omega)}ds+\frac{rm^*}{|\Omega|}\nonumber\\
&\le \bar{C}, ~t\in(0,T_{\max})
\end{align}
with some $\bar{C}>0$. This concludes $T_{\max}=\infty$ by Lemma \ref{lemma1}, i.e., the classical solution $(u,v)$ is globally bounded.  \qquad\medskip$\Box$
\section{Regularization problem}
To deal with the global existence-boundedness of generalized solution to \eqref{p} for $k>1$, we introduce an appropriate regularization problem related to \eqref{p}
\begin{eqnarray}\label{pp}
 \left\{
\begin{array}{llll}
u_{\epsilon t}=\Delta u_\epsilon-\chi\nabla\cdot(\frac{u_\epsilon}{v_\epsilon}\nabla v_\epsilon)
+ru_\epsilon-\mu u_\epsilon^k-\epsilon u_\epsilon^{k+1}, & x\in\Omega,~~t>0,\\[4pt]
\displaystyle 0=\Delta v_\epsilon-v_\epsilon+u_\epsilon,& x\in\Omega,~~t>0,\\[4pt]
  \displaystyle \frac{\partial u_\epsilon}{\partial {\nu}}=\frac{\partial v_\epsilon}
  {\partial {\nu}}=0 ,& x\in\partial\Omega,~~t>0,\\[4pt]
  \displaystyle u_\epsilon(x,0)=u_0(x),  &x\in\Omega
\end{array}\right. \end{eqnarray}
with $\epsilon\in(0,1)$.
The local classical solution of the regularization problem \eqref{pp} with general $k>1$ can be obtained in the similar arguments \cite{FF}. That is:
\begin{lemma}\label{lemma400}
Assume that $u_0$ satisfies \eqref{con}. Let $k>1$, $r,\chi,\mu>0$. Then
for each $\epsilon\in(0,1)$ there exist $T_{\max,\epsilon}\in(0,+\infty]$ and a unique pair $(u_\epsilon,v_\epsilon)$ of functions
 \begin{equation*}
   \left\{
    \begin{aligned}
    &u_\epsilon\in C^0(\overline\Omega\times[0,T_{\max,\epsilon}))\cap C^{2,1}
     (\overline{\Omega}\times(0,T_{\max,\epsilon})),\\[4pt]
    &v_\epsilon\in  C^{2,0}
     (\overline{\Omega}\times(0,T_{\max,\epsilon})),
    \end{aligned}
    \right.
\end{equation*}
satisfying \eqref{p} in the classical sense with $u_\epsilon,v_\epsilon>0$ in $\overline\Omega\times(0,T_{\max,\epsilon})$. Moreover, either $ T_{\max,\epsilon}=\infty$, or ${\limsup}_{t\rightarrow T_{\max,\epsilon}}\| u_\epsilon(\cdot,t)\|_{L^\infty(\Omega)}=\infty$, or ${\liminf}_{t\rightarrow T_{\max,\epsilon}} \inf_{x\in\Omega} v_\epsilon(x,t)=0$.   $\Box$
\end{lemma}

Let $(u_\epsilon,v_\epsilon)$ is the local classical solution to system \eqref{pp} for each $\epsilon\in(0,1)$ and $k>1$. Without loss of generality, assume that $T_{\max,\epsilon}>1$ for each $\epsilon\in(0,1)$. Then we have the following a priori estimates.

\begin{lemma}\label{lemma33}
With $k>1$ and $\mu,r,\chi>0$, it holds for each $\epsilon>0$ that
\begin{align}\label{lem331}
&\int_\Omega u_\epsilon dx\le m^*,~~t\in(0,T_{\max,\epsilon}),\\
\label{lem332}
&\int_0^T\int_\Omega u_\epsilon^kdxds\le M_1(1+T),~~T\in(0,T_{\max,\epsilon}-1)\\
&\int_t^{t+1}\int_\Omega u_\epsilon^kdxds\le M_1,~~t\in(0,T_{\max,\epsilon})
\label{lem333}
\end{align}
with $m^*,M_1>0$ defined in Lemma \ref{lemma2}.
\end{lemma}
{\bf Proof.}\ Integrate \eqref{pp}$_1$ over $\Omega$ with the H\"{o}lder inequality to know that
\begin{align}
\frac{d}{dt}\int_\Omega u_\epsilon dx& =r\int_\Omega u_\epsilon dx-\mu\int_\Omega
u_\epsilon^kdx-\epsilon\int_\Omega
u_\epsilon^{k+1}dx\label{331}\\
&\le r\int_\Omega u_\epsilon dx-\frac{\mu}{|\Omega|^{k-1}}\big(\int_\Omega
u_\epsilon dx\big)^k,~~t>0.\label{332}
\end{align}
We get \eqref{lem331} by the Bernoulli inequality with \eqref{332}.
The estimates \eqref{lem332} and \eqref{lem333} come from by integrating \eqref{331} with \eqref{lem331}.
  \qquad $\Box$\medskip

In order to deal with the global existence of classical solution to \eqref{pp} for each $\epsilon\in(0,1)$, we give a uniform-in-time lower bound of $v_\epsilon$  for all $(x,t)\in\Omega\times(0,T_{\max,\epsilon})$.
 \begin{lemma}\label{lemma40}
Let $k>1$, $\mu>0$ and $r,\chi>0$ satisfy \eqref{lem51}. Then for each $\epsilon\in(0,1)$ there exists some $\delta_1>0$ such that
\begin{align}\label{lem401}
v_\epsilon(x,t)\ge\delta_1~~{\rm for~ all}~(x,t)\in\Omega\times(0,T_{\max,\epsilon}).
\end{align}
\end{lemma}
{\bf Proof.}\
Since $u_\epsilon\in C^0(\bar{\Omega}\times[0,T_{\max,\epsilon}))$ by Lemma \ref{lemma400} for each $\epsilon\in(0,1)$, we know by a continuous argument that there exists $t_0\in(0,1)$  such that
 \begin{align}\label{400}
 \int_\Omega u_\epsilon(x,t)dx\ge \frac{1}{2}\int_\Omega u_0dx,~~t\in(0,t_0].
 \end{align}

If $m>0$, a direct computation with \eqref{pp}$_1$ shows
\begin{align}\label{401}
\frac{1}{m}\frac{d}{dt}\int_\Omega u_\epsilon^{-m}dx
&=-\int_\Omega u_\epsilon^{-m-1}[\Delta u_\epsilon-\chi\nabla\cdot(\frac{u_\epsilon}{v_\epsilon}\nabla v_\epsilon)]dx\nonumber\\
&~~~-r\int_\Omega u_\epsilon^{-m}dx+\mu\int_\Omega u_\epsilon^{-m-1+k}dx +\epsilon\int_\Omega u_\epsilon^{-m+k}dx\nonumber\\
&=-(m+1)\int_\Omega u_\epsilon^{-m-2}|\nabla u_\epsilon|^2dx+\chi(m+1)\int_\Omega \frac{u_\epsilon^{-m-1}}{v_\epsilon}\nabla u_\epsilon\cdot\nabla v_\epsilon dx \nonumber\\
&~~~-r\int_\Omega u_\epsilon^{-m}dx+\mu\int_\Omega u_\epsilon^{-m-1+k}dx +\epsilon\int_\Omega u_\epsilon^{-m+k}dx,~~t\in(t_0,T_{\max,\epsilon})
\end{align}
for all $\epsilon\in(0,1)$. If $m=\frac{4a}{(\chi-a)^2}$ with $a\in(0,\chi)$ and $r,\chi>0$ satisfying \eqref{con},  there exists some $c_0>0$ such that $-(r-\frac{(m+1)a}{m})\le-c_0<0$ by a similar argument as that in Lemma \ref{lemma5}, and hence it holds from \eqref{401} that
\begin{align}\label{404}
\frac{1}{m}\frac{d}{dt}\int_\Omega u_\epsilon^{-m}dx\le -c_0\int_\Omega u_\epsilon^{-m}dx+\mu\int_\Omega u_\epsilon^{-m-1+k}dx +\int_\Omega u_\epsilon^{-m+k}dx,~t\in(t_0,T_{\max,\epsilon})
\end{align}
for all  $\epsilon\in(0,1)$.

If $k-m\le 0$, we have from \eqref{404} by the Young's inequality that
\begin{align}\label{405}
\frac{1}{m}\frac{d}{dt}\int_\Omega u_\epsilon^{-m}dx&\le -\frac{c_0}{2}\int_\Omega u_\epsilon^{-m}dx+\Big(1+(\frac{4}{c_0})^\frac{m-k}{k}+\mu(\frac{4\mu}{c_0})^\frac{m+1-k}{k-1}\Big)|\Omega|.
\end{align}

Similarly, if $k-1-m\le 0<k-m$, we get
\begin{align}\label{406}
\frac{1}{m}\frac{d}{dt}\int_\Omega u_\epsilon^{-m}dx&\le -\frac{c_0}{2}\int_\Omega u_\epsilon^{-m}dx+\int_\Omega u_\epsilon ^kdx+\Big(2+\mu(\frac{2\mu}{c_0})^\frac{m+1-k}{k-1}\Big)|\Omega|,
\end{align}
and if $k-1-m> 0$,
\begin{align}\label{407}
\frac{1}{m}\frac{d}{dt}\int_\Omega u_\epsilon^{-m}dx&\le -c_0\int_\Omega u_\epsilon^{-m}dx+\int_\Omega u_\epsilon ^kdx+\Big(\mu(2\mu)^\frac{k-1-m}{m+1}+2^\frac{k-m}{m}\Big)|\Omega|.
\end{align}
The estimates \eqref{405}--\eqref{407} yield for $k>1$, $\mu>0$ and $r,\chi>0$ satisfying \eqref{lem51} that
\begin{align}\label{408}
\frac{1}{m}\frac{d}{dt}\int_\Omega u_\epsilon^{-m}dx&\le -\frac{c_0}{2}\int_\Omega u_\epsilon^{-m}dx+\int_\Omega u_\epsilon ^kdx+C_{11},~t\in(t_0,T_{\max,\epsilon})
\end{align}
with some $C_{11}>0$ for all $\epsilon\in(0,1)$. By a similar discussion as that in Lemma \ref{lemma5} with \eqref{lem333}, we conclude the uniformly lower bound estimate \eqref{lem401} with some $\delta_1>0$. \qquad$\Box$\medskip

If $k>2-\frac{2}{n}$, the local classical solution $(u_\epsilon,v_\epsilon)$ for each $\epsilon\in(0,1)$ is in fact global.
\begin{lemma}\label{lemma30}
Let $k>2-\frac{2}{n}$, $\mu>0$ and $r,\chi>0$ satisfy \eqref{lem51}. Then for each $\epsilon\in(0,1)$ the system \eqref{pp} possesses a global classical solution $(u_\epsilon,v_\epsilon)$.
\end{lemma}
{\bf Proof.}\
For each $\epsilon\in(0,1)$, let $(u_\epsilon, v_\epsilon)$ be the local classical solution  to the regularization problem \eqref{pp} with general $k>1$. Since $k+1>\frac{3n-2}{n}$ for $k>2-\frac{2}{n}$, replacing $\delta_0$ by $\delta_1$ in \eqref{251} and \eqref{tt3}, we can prove for each $\epsilon\in(0,1)$ that the solution $(u_\epsilon,v_\epsilon)$ is global by a similar arguments in Lemma \ref{lemma25} and \eqref{tt1}--\eqref{tt3}.
 \qquad\medskip$\Box$

\begin{corollary}\label{re4}
Let $k>2-\frac{2}{n}$, $\mu>0$ and $r,\chi>0$ satisfy \eqref{lem51}. Then the estimates in Lemmas \ref{lemma33} and \ref{lemma40} are  valid with $T_{\max,\epsilon}=\infty$ and for all $\epsilon\in(0,1)$. For convenience, we omits the new marks on these esitimates.
\end{corollary}
{\bf Proof.}\
For each $\epsilon\in(0,1)$, it is shown from Lemma \ref{lemma30} that the regularization problem \eqref{pp} possesses a global classical $(u_\epsilon,v_\epsilon)$. This yields $T_{\max,\epsilon}=\infty$ for all $\epsilon\in(0,1)$. Since the constants $m^*, M_1>0$ are not dependent on $\epsilon$, the estimates \eqref{lem331}--\eqref{lem333} in Lemma \ref{lemma33} are valid for all $\epsilon\in(0,1)$. In addition, we know $u_\epsilon\in C^0(\bar{\Omega}\times[0,1))$ for all $\epsilon\in(0,1)$, which concludes the estimate \eqref{400} with some $t_0>0$ independent of $\epsilon$. A similar arguments (the constants there are all independent of $\epsilon$) from \eqref{401} to \eqref{408}, and in Lemma \ref{lemma5} with \eqref{lem333} indicate the estimate \eqref{lem401} for all $\epsilon\in(0,1)$.\qquad $\Box$\medskip

Now, we deal with a spatio-temporal integral estimate on $ \nabla \ln v_\epsilon$ for all $\epsilon\in(0,1)$.
\begin{lemma}\label{lemma41}
Let $k>2-\frac{2}{n}$, $\mu>0$ and  $r, \chi>0$ satisfy \eqref{lem51}. Then
\begin{align}\label{lem411}
\int_0^T\int_\Omega \frac{|\nabla v_\epsilon|^2}{v_\epsilon^2}dxds\le |\Omega|T,~~T>0
\end{align}
for all $\epsilon\in(0,1)$.
\end{lemma}
{\bf Proof.}\
Multiply \eqref{pp}$_2$ by $\frac{1}{v_\epsilon}$ and integrate by part to get
\begin{align}\label{411}
0&=\int_\Omega \frac{1}{v_\epsilon}[\Delta v_\epsilon-v_\epsilon+u_\epsilon]dx=\int_\Omega\frac{|\nabla v_\epsilon|^2}{v_\epsilon^2}dx-|\Omega|+\int_\Omega \frac{u_\epsilon}{v_\epsilon}dx,~t>0.
\end{align}
 This yields conclusion \eqref{lem411} by integrating \eqref{411} from $0$ to $T$ for all $\epsilon\in(0,1)$.  \qquad$\Box$\medskip

We proceed to derive another spatio-temporal integral estimate on $u_\epsilon$ for all $\epsilon\in(0,1)$.

\begin{lemma}\label{lemma42}
For $k>2-\frac{2}{n}$, $\mu>0$ and $r,\chi>0$ satisfying \eqref{lem51}, there exists some $M_3>0$ such that
\begin{align}\label{lem421}
\int_0^T\int_\Omega \frac{|\nabla u_\epsilon|^2}{(1+u_\epsilon)^2}dxdt\le M_3(1+T),~~T>0
\end{align}
for all $\epsilon\in(0,1)$.
\end{lemma}
{\bf Proof.}\ A direct calculation with \eqref{pp}$_1$ shows
\begin{align}\label{421}
\frac{d}{dt}\int_\Omega \ln(1+u_\epsilon) dx&=\int_\Omega \frac{1}{1+u_\epsilon}[\Delta u_\epsilon-\chi\nabla\cdot(\frac{u_\epsilon}{v_\epsilon}\nabla v_\epsilon)+ru_\epsilon-\mu u_\epsilon^k-\epsilon u_\epsilon^{k+1}]dx\nonumber\\
&=\int_\Omega \frac{|\nabla u_\epsilon|^2}{(1+u_\epsilon)^2}dx-\chi\int_\Omega \frac{u_\epsilon}{(1+u_\epsilon)^2v_\epsilon}\nabla u_\epsilon\cdot v_\epsilon dx\nonumber\\
&~~~+r\int_\Omega \frac{u_\epsilon}{1+u_\epsilon} dx-\mu\int_\Omega \frac{u_\epsilon^k}{1+u_\epsilon}dx-\epsilon \int_\Omega\frac{u_\epsilon^{k+1}}{1+u_\epsilon}dx,~~t>0.
\end{align}
By Young's inequality, we have
\begin{align}\label{422}
\chi\int_\Omega \frac{u_\epsilon}{(1+u_\epsilon)^2v_\epsilon}\nabla u_\epsilon\cdot v_\epsilon dx\le\frac{1}{2}\int_\Omega \frac{|\nabla u_\epsilon|^2}{(1+u_\epsilon)^2}dx+\frac{\chi^2}{2}\int_\Omega \frac{|\nabla v_\epsilon|^2}{v_\epsilon^2}dx.
\end{align}
It is known from \eqref{421} and \eqref{422} that
\begin{align}\label{423}
\int_\Omega \frac{|\nabla u_\epsilon|^2}{(1+u_\epsilon)^2}dx&\le 2\frac{d}{dt}\int_\Omega \ln(1+u_\epsilon)dx+\chi^2\int_\Omega \frac{|\nabla v_\epsilon|^2}{v_\epsilon^2}dx\nonumber\\
&~~~-2r\int_\Omega \frac{u_\epsilon}{1+u_\epsilon} d+2\mu\int_\Omega \frac{u_\epsilon^k}{1+u_\epsilon}dx+2\epsilon \int_\Omega\frac{u_\epsilon^{k+1}}{1+u_\epsilon}dx,~~t>0.
\end{align}
Combining \eqref{423} with \eqref{lem331}, \eqref{lem332} and \eqref{lem411}, then we get with the fact $0<\ln(1+a)\le a$ for $a>0$ that
\begin{align}\label{424}
\int_0^T\int_\Omega \frac{|\nabla u_\epsilon|^2}{(1+u_\epsilon)^2}dxds&\le 2\int_\Omega \ln (1+u_\epsilon(\cdot,t))dx-2\int_\Omega \ln(1+u_0)dx+\chi^2\int_0^T\int_\Omega \frac{|\nabla v_\epsilon|^2}{v_\epsilon^2}dxds\nonumber\\
&~~-2r\int_0^T\int_\Omega \frac{u_\epsilon}{1+u_\epsilon} dxds+\mu\int_0^T\int_\Omega \frac{u_\epsilon^{k}}{1+u_\epsilon}dxds+\epsilon\int_0^T\int_\Omega \frac{u_\epsilon^{k+1}}{1+u_\epsilon}dxds\nonumber\\
&\le 2\int_\Omega u_\epsilon(\cdot,t) dx+(1+\mu)\int_0^T\int_\Omega u_\epsilon^kdxds+\chi^2\int_0^T\int_\Omega \frac{|\nabla v_\epsilon|^2}{v_\epsilon^2}dxds\nonumber\\
&\le M_3(1+T),~~T>0
\end{align}
with $M_{3}>0$.  \qquad$\Box$\medskip

Next, we deal with the estimate on the time derivative of $\ln(1+u_\epsilon)$ for all $\epsilon\in(0,1)$.
\begin{lemma}\label{lemma43}
For $k>2-\frac{2}{n}$, $\mu>0$ and $r,\chi>0$ satisfying \eqref{lem51},  there exists $M_4>0$ such that
\begin{align}\label{lem431}
\int_0^T\Big\|\frac{d}{dt}\ln(1+u_\epsilon)\Big\|_{(W_0^{n+1,2}(\Omega))^*}ds\le M_4(1+T),~~T>0
\end{align}
for all $\epsilon\in(0,1)$.
\end{lemma}
{\bf Proof.}\ Let $\phi\in W_0^{n+1,2}(\Omega)$. Then we have from \eqref{pp}$_1$ that
\begin{align}\label{431}
\int_\Omega \frac{d}{dt}&\ln(1+u_\epsilon)\phi dx=\int_\Omega \frac{1}{(1+u_\epsilon)}\phi[\Delta u_\epsilon-\chi\nabla\cdot(\frac{u_\epsilon}{v_\epsilon}\nabla v_\epsilon)+r u_\epsilon-\mu u_\epsilon^k-\epsilon u_\epsilon^{k+1}]dx\nonumber\\
&=-\int_\Omega\nabla(\frac{\phi}{1+u_\epsilon})\cdot(\nabla u_\epsilon-\chi\frac{u_\epsilon}{v_\epsilon}\nabla v_\epsilon)dx+r\int_\Omega \frac{u_\epsilon}{1+u_\epsilon}\phi dx-\mu\int_\Omega \frac{u_\epsilon^k}{1+u_\epsilon}dx-\epsilon\int_\Omega \frac{u_\epsilon^{k+1}}{1+u_\epsilon}dx\nonumber\\
&=\int_\Omega \frac{|\nabla u_\epsilon|^2}{(1+u_\epsilon)^2}\phi dx-\int_\Omega\frac{\nabla u_\epsilon\cdot\nabla \phi}{1+u_\epsilon}dx-\chi\int_\Omega\frac{u_\epsilon}{(1+u_\epsilon)^2v_\epsilon}\nabla u_\epsilon\cdot\nabla v_\epsilon \phi dx\nonumber\\
&~~+\chi\int_\Omega\frac{u_\epsilon}{(1+u_\epsilon)v_\epsilon}\nabla v_\epsilon\cdot\nabla \phi dx+r\int_\Omega \frac{u_\epsilon}{1+u_\epsilon}\phi dx-\mu\int_\Omega \frac{u_\epsilon^k}{1+u_\epsilon}\phi dx-\epsilon\int_\Omega \frac{u_\epsilon^{k+1}}{1+u_\epsilon}\phi dx\nonumber\\
&\le \big(\int_\Omega \frac{|\nabla u_\epsilon|^2}{(1+u_\epsilon)^2}dx\big)\|\phi\|_{L^\infty(\Omega)}+\big(\int_\Omega \frac{|\nabla u_\epsilon|^2}{(1+u_\epsilon)^2}dx\big)^\frac{1}{2}\|\nabla\phi\|_{L^2(\Omega)}\nonumber\\
&~~+{\chi}\big(\int_\Omega \frac{|\nabla v_\epsilon|^2}{v_\epsilon^2} dx\big)^\frac{1}{2}\|\nabla \phi\|_{L^2(\Omega)}+\big(\int_\Omega \frac{|\nabla u_\epsilon|^2}{(1+u_\epsilon)^2}dx\big)^\frac{1}{2}\big(\int_\Omega \frac{|\nabla v_\epsilon|^2}{v_\epsilon^2}dx\big)^\frac{1}{2}\|\phi\|_{L^\infty(\Omega)}\nonumber\\
&~~+\big(r+(1+\mu)\int_\Omega u_\epsilon^kdx\big)\|\phi\|_{L^\infty(\Omega)},~t>0
\end{align}
by the H\"{o}lder inequality. Since $W_0^{n+1,2}(\Omega)\hookrightarrow W^{1,\infty}(\Omega)$, it is known by Young's inequality with \eqref{431}  that
\begin{align}\label{432}
\Big|\int_\Omega \frac{d}{dt}\ln(1+u_\epsilon)\phi dx\Big|&\le C_{12}\Big(1+\int_\Omega u_\epsilon^kdx+\int_\Omega \frac{|\nabla u_\epsilon|^2}{(1+u_\epsilon)^2}dx+\int_\Omega \frac{|\nabla v_\epsilon|^2}{v_\epsilon^2}dx\Big)\|\phi\|_{W_0^{n+1,2}(\Omega)}
\end{align}
with $C_{12}>0$ for $t>0$.
Integrating \eqref{432} from $0$ to $T$,  we obtain from \eqref{lem332}, \eqref{lem411}, \eqref{lem421} and \eqref{432} that
\begin{align*}
\int_0^T\Big\|\frac{d}{dt}\ln(1+u_\epsilon&)\Big\|_{(W_0^{n+1,2}(\Omega))^*}ds\le \sup_{\phi\in W_0^{n+1,2}(\Omega),\|\phi\|_{W_0^{n+1,2}(\Omega)}\le 1}\int_0^T\Big|\int_\Omega \frac{d}{dt}\ln(1+u_\epsilon)\phi dx\Big|ds\\
&\le C_{12}\Big(T+\int_0^T\int_\Omega u_\epsilon^kdxds+\int_0^T\int_\Omega \frac{|\nabla u_\epsilon|^2}{(1+u_\epsilon)^2}dxds+\int_0^T\int_\Omega \frac{|\nabla v_\epsilon|^2}{v_\epsilon^2}dxds \Big)\\
&\le M_4(1+T),~~T>0
\end{align*}
with  some  $M_4>0$. The proof is complete.
\qquad$\Box$\medskip

Based on Lemma \ref{lemma40}, we further prove the following estimates on $v_\epsilon$ for all $\epsilon\in(0,1)$.
\begin{lemma}\label{lemma44}
Let $k>2-\frac{2}{n}$, $\mu>0$ and $r,\chi>0$ satisfy \eqref{lem51}. Then for $q\in(2,\frac{nk}{n-1})$ there exists $M_5>0$ such that
\begin{align}\label{lem442}
\int_0^T\int_\Omega\big|\frac{\nabla v_\epsilon}{v_\epsilon}\big|^qdxds\le M_5(1+T),~T>0
\end{align}
for all $\epsilon\in(0,1)$.
\end{lemma}
{\bf Proof.}\
Similar argument for \eqref{252} and \eqref{2521} in Lemma \ref{lemma25}, we know for $r\in(1,\frac{n}{n-1})$ that
\begin{align}\label{442}
\|v_\epsilon(\cdot,t)\|_{W^{1,r}(\Omega)}\le C_{BS}\|\Delta v_\epsilon\|_{L^1(\Omega)}\le C_{BS}\|v_{\epsilon}-u_\epsilon\|_{L^1(\Omega)}\le 2C_{BS}m^*,~t>0
\end{align}
and for $m\ge 1$ that
\begin{align}\label{443}
 \|v_\epsilon(\cdot,t)\|_{W^{2,m}(\Omega)}\le C_5\|u_\epsilon(\cdot,t)\|_{L^{m}(\Omega)},~t>0
\end{align}
for all $\epsilon\in(0,1)$. Let $q\in(\frac{(n+1)k}{n},\frac{nk}{n-1})$. Then $r:=\frac{n(q-k)}{k}\in(1,\frac{n}{n-1})$, and hence by the Gagliardo-Nirenberg inequality with \eqref{443}, we know
\begin{align*}
\|\nabla v_\epsilon(\cdot,t)\|_{L^{q}(\Omega)}&\le C_{GN}\| v_\epsilon(\cdot,t)\|_{W^{2,k}(\Omega)}^a\|\nabla v_\epsilon(\cdot,t)\|_{W^{1,r}(\Omega)}^{1-a}\nonumber\\
&\le C_{5}C_{GN}\| u_\epsilon(\cdot,t)\|_{L^k(\Omega)}^a\|\nabla v_\epsilon(\cdot,t)\|_{W^{1,r}(\Omega)}^{1-a},~t>0
\end{align*}
with some $C_{GN}>0$ for all $\epsilon\in(0,1)$, where $a=\frac{\frac{n}{r}-\frac{n}{q}}{1-\frac{n}{q}+\frac{n}{r}}\equiv\frac{k}{q}\in(0,1)$.
This together with \eqref{442} and \eqref{lem332} indicates
\begin{align}\label{444}
\int_0^T\int_\Omega|{\nabla v_\epsilon}|^qdxdt\le C_{13}(1+T),~T>0
\end{align}
with some $C_{13}>0$. The proof is complete by \eqref{444} and \eqref{lem401} with $M_5=\frac{C_{13}}{\delta_1^q}>0$.
\qquad$\Box$\medskip

We now perform a subsequence extraction procedure to obtain a limit object $(u,v)$, i.e.,  a generalized solution to the problem \eqref{p}.
\begin{lemma}\label{lemma45}
Let $k>2-\frac{2}{n}$, $\mu>0$ and $r,\chi>0$ satisfy \eqref{lem51}. Then for $p\in(1,k)$ with $q\in(2,\frac{nk}{k-1})$ there exist $u\in L_{\rm loc}^1(\Omega\times(0,\infty))$ and $v\in L_{\rm loc}^1((0,\infty),W^{1,1}(\Omega))$ such that
\begin{align}\label{lem451}
\ln(1+u_\epsilon)&\rightharpoonup \ln(1+u),  &in ~L_{\rm loc}^2([0,\infty);W^{1,2}(\Omega)),\\
\label{lem452}
u_\epsilon&\rightharpoonup u,   &~in~L_{\rm loc}^k(\Omega\times[0,\infty)),\\
\label{lem453}
u_\epsilon&\rightarrow u,       &a.e. ~in ~\Omega\times(0,\infty)~and ~L_{\rm loc}^p(\Omega\times[0,\infty))\\
\label{lem454}
v_\epsilon&\rightarrow v,     &a.e. ~in ~\Omega\times(0,\infty)~and~in~L_{\rm loc}^1([0,\infty);W^{1,1}(\Omega))\\
\label{lem455}
v_\epsilon&\rightharpoonup v, &~in~L_{\rm loc}^q([0,\infty);W^{1,q}(\Omega))\\
\label{lem456}
\frac{ u_\epsilon^2|\nabla v_\epsilon|^2}{(1+u_\epsilon)^2v_\epsilon^2}&\rightarrow \frac{u^2|\nabla v|^2}{(1+u)^2v^2},&~in~L_{\rm loc}^1(\Omega\times[0,\infty))
\end{align}
for $\epsilon=\epsilon_j\searrow 0$.
\end{lemma}
{\bf Proof.}\
Let $T>0$. The conclusions \eqref{lem451}, \eqref{lem452} and \eqref{lem455} are the direct results from \eqref{lem421}, \eqref{lem332} and \eqref{444}.
Since $W^{1,2}(\Omega)\hookrightarrow\hookrightarrow L^2(\Omega)$, we have by the Aubin-Lions lemma with \eqref{lem421} and \eqref{lem431} that $\ln(1+u_\epsilon)\rightarrow \ln(1+u)$ in $L^2(\Omega\times(0,T))$, and moreover $u_\epsilon\rightarrow u$ a.e. in $\Omega\times(0,T)$, as  $\epsilon=\epsilon_j\searrow 0$. For $p\in(1,k)$, it is known from \eqref{lem332} that
\begin{align*}
\int_0^T\int_\Omega u_\epsilon^pdxdt\le |\Omega|T +\int_0^T\int_\Omega u_\epsilon^kdxdt\le C(1+T),~T>0
\end{align*}
by Young's inequality for all $\epsilon\in(0,1)$, i.e., ${\{u_\epsilon^p\}_{\epsilon\in(0,1)}}\subset L_{\rm loc}^\frac{k}{p}(\Omega\times[0,\infty))$. This together with $u_\epsilon\rightarrow u$ a.e. in $\Omega\times (0,T)$ indicates \eqref{lem453} by the Vitali convergence theorem.
 The estimates \eqref{443} and \eqref{lem453} imply that there exists some nonnegative $v$ defined on $\Omega\times (0,T)$ such that \eqref{lem454} holds. Consequently, we note from \eqref{lem442} for $q\in(2,\frac{nk}{k-1})$ that $
\frac{ |\nabla v_\epsilon|^2}{v_\epsilon^2}\rightharpoonup \frac{|\nabla v|^2}{v^2}$ in $L_{\rm loc}^\frac{q}{2}(\Omega\times[0,\infty))$, which concludes \eqref{lem456} due to $\frac{u_\epsilon}{1+u_\epsilon}\rightarrow \frac{u}{1+u}$ in $L_{\rm loc}^m(\Omega\times[0,\infty))$ for every $m>1$ by \eqref{lem453} as $\epsilon=\epsilon_j\searrow 0$.  \qquad$\Box$\medskip

\section{Global existence and boudedness to generalized solution}
In this section we begin with proving that the function $(u,v)$ determined in Lemma \ref{lemma45} just is the global generalized solution of \eqref{p}.\\[6pt]
{\bf Proof the Theorem \ref{th2}.}\
For $k> 2-\frac{1}{n}$, we will firstly demonstrate that the function $(u,v)$ obtained in Lemma  \ref{lemma45} is a very weak subsolution of \eqref{p} in $\Omega\times(0,T)$ for $T>0$. Let $\varphi$ satisfy \eqref{def13}. Multiplying  \eqref{pp}$_1$ by $\varphi$ and integrating by parts, then we have for all $\epsilon\in(0,1)$ that
\begin{align}\label{t11}
-\int_0^T\int_\Omega u_\epsilon\varphi_t -\int_\Omega u_0\varphi(\cdot,0)&= \int_0^T\int_\Omega u_\epsilon\Delta \varphi+\chi\int_0^T\int_\Omega \frac{u_\epsilon}{v_\epsilon}\nabla v_\epsilon\cdot\nabla\varphi\nonumber\\
&~+r\int_0^T\int_\Omega u_\epsilon\varphi-\mu\int_0^T\int_\Omega u_\epsilon^k\varphi-\epsilon\int_0^T\int_\Omega u_\epsilon^{k+1}\varphi.
\end{align}
By \eqref{lem452}, we know
\begin{align}\label{t12}
-\int_0^T\int_\Omega u_\epsilon\varphi_t&\rightarrow -\int_0^T\int_\Omega u\varphi_t,\\
\label{t13}
\int_0^T\int_\Omega u_\epsilon\Delta \varphi&\rightarrow \int_0^T\int_\Omega u\Delta \varphi,\\
\label{t14}
r\int_0^T\int_\Omega u_\epsilon\varphi&\rightarrow r\int_0^T\int_\Omega u\varphi
\end{align}
as $\epsilon=\epsilon_j \searrow 0$.
Since $\frac{nk}{n-1}>\frac{k}{k-1}$ for $k>2-\frac{1}{n}$, we know by \eqref{lem453} with \eqref{lem455} that
\begin{align}\label{t16}
\chi\int_0^T\int_\Omega u_\epsilon\frac{\nabla v_\epsilon}{v_\epsilon}\cdot\nabla\varphi\rightarrow \chi\int_0^T\int_\Omega u\frac{\nabla v}{v}\cdot\nabla\varphi~~{\rm as}~ \epsilon=\epsilon_j \searrow 0.
\end{align}
Consequently, in view of \eqref{t12}--\eqref{t16} with the Fatou lemma and the positivity of $\epsilon\int_0^T\int_\Omega u_\epsilon^{k+1}\varphi$ for $\epsilon\in(0,1)$, we obtain
\begin{align}\label{t17}
\mu\int_0^T\int_\Omega u^k\varphi&\le \mu\liminf_{\epsilon=\epsilon_j\searrow 0}\int_0^T\int_\Omega u_\epsilon^k\varphi\nonumber\\
&=\int_0^T\int_\Omega u\varphi_t +\int_\Omega u_0\varphi(\cdot,0)+\int_0^T\int_\Omega u\Delta \varphi\nonumber\\&~~+\chi\int_0^T\int_\Omega u\frac{\nabla v}{v}\cdot\nabla\varphi+r\int_0^T\int_\Omega u\varphi.
\end{align}

Take $\psi$ satisfying \eqref{def14}. Multiply \eqref{pp}$_2$ by $\psi$ and integrate by parts, then we
\begin{align}\label{t18}
-\int_0^T\int_\Omega v_\epsilon\psi_t-\int_\Omega v_0\psi(\cdot,0)+\int_0^T\int_\Omega \nabla v_\epsilon\cdot\nabla \psi+\int_0^T\int_\Omega v_\epsilon\psi=\int_0^T\int_\Omega u_\epsilon\psi.
\end{align}
According to \eqref{lem454} and \eqref{lem452}, we get \eqref{def12} by taking $\epsilon=\epsilon_j \searrow 0$. This together with \eqref{t17} indicates that $(u,v)$ is a very weak subsolution of \eqref{p}.

Taking $\varphi$ in \eqref{def13} and multiplying \eqref{p}$_1$ by $\frac{\varphi}{1+u_\epsilon}$, we have
\begin{align}\label{t19}
-\int_0^T\int_\Omega \ln(1+u_\epsilon)\varphi_t&-\int_\Omega \ln(1+u_0)\varphi(\cdot,0)=\int_0^T\int_\Omega \frac{|\nabla u_\epsilon|^2}{(1+u_\epsilon)^2}\varphi-\chi\int_0^T\int_\Omega \frac{u_\epsilon}{(1+u_\epsilon)^2v_\epsilon}\nabla u_\epsilon\cdot\nabla v_\epsilon \varphi\nonumber\\
&~~-\int_0^T\int_\Omega\frac{\nabla u_\epsilon \cdot\nabla \varphi}{1+ u_\epsilon}+\chi\int_0^T\int_\Omega\frac{u_\epsilon}{(1+u_\epsilon)v_\epsilon} \nabla v_\epsilon\cdot\nabla\varphi\nonumber\\
&~~
+r\int_0^T\int_\Omega \frac{u_\epsilon}{1+u_\epsilon}\varphi-\mu\int_0^T\int_\Omega \frac{u_\epsilon^k}{1+u_\epsilon}\varphi-\epsilon\int_0^T\int_\Omega \frac{u_\epsilon^{k+1}}{1+u_\epsilon}\varphi\nonumber\\
&=\int_0^T\int_\Omega \Big(\frac{\nabla u_\epsilon}{1+u_\epsilon}-\frac{\chi u_\epsilon\nabla v_\epsilon}{2(1+u_\epsilon)v_\epsilon}\Big)^2\varphi-\frac{\chi^2}{4}\int_0^T\int_\Omega\frac{u_\epsilon^2|\nabla v_\epsilon|^2}{(1+u_\epsilon)^2v_\epsilon^2}\varphi\nonumber\\
&~~-\int_0^T\int_\Omega\frac{\nabla u_\epsilon \cdot\nabla \varphi}{1+ u_\epsilon}+\chi\int_0^T\int_\Omega\frac{u_\epsilon}{(1+u_\epsilon)v_\epsilon} \nabla v_\epsilon\cdot\nabla\varphi\nonumber\\
&~~
+r\int_0^T\int_\Omega \frac{u_\epsilon}{1+u_\epsilon}\varphi-\mu\int_0^T\int_\Omega \frac{u_\epsilon^k}{1+u_\epsilon}\varphi-\epsilon\int_0^T\int_\Omega \frac{u_\epsilon^{k+1}}{1+u_\epsilon}\varphi.
\end{align}
By \eqref{lem452}, we know
\begin{align}\label{t110}
-\int_0^T\int_\Omega \ln(1+u_\epsilon)\varphi_t&\rightarrow -\int_0^T\int_\Omega \ln(1+u)\varphi_t,\\
\label{t111}
r\int_0^T\int_\Omega \frac{u_\epsilon}{1+u_\epsilon}\varphi&\rightarrow r\int_0^T\int_\Omega \frac{u}{1+u}\varphi,\\
\label{t112}
-\mu\int_0^T\int_\Omega \frac{u_\epsilon^k}{1+u_\epsilon}\varphi&\rightarrow -\mu\int_0^T\int_\Omega \frac{u^k}{1+u}\varphi
\end{align}
as $\epsilon=\epsilon_j\searrow 0$, whereas \eqref{lem451} implies that
\begin{align}\label{t113}
-\int_0^T\int_\Omega\frac{\nabla u_\epsilon \cdot\nabla \varphi}{1+ u_\epsilon}\rightarrow -\int_0^T\int_\Omega\frac{\nabla u \cdot\nabla \varphi}{1+ u}
\end{align}
as $\epsilon=\epsilon_j\searrow 0$. It follows from  \eqref{lem456} that
\begin{align}\label{t114}
-\frac{\chi^2}{4}\int_0^T\int_\Omega\frac{u_\epsilon^2|\nabla v_\epsilon|^2}{(1+u_\epsilon)^2v_\epsilon^2}\varphi&\rightarrow-\frac{\chi^2}{4}\int_0^T\int_\Omega\frac{u^2|\nabla v|^2}{(1+u)^2v^2}\varphi,\\
\label{t115}
\chi\int_0^T\int_\Omega\frac{u_\epsilon}{(1+u_\epsilon)v_\epsilon}\nabla v_\epsilon \cdot\nabla\varphi&\rightarrow \chi\int_0^T\int_\Omega\frac{u}{(1+u)v} \nabla v\cdot\nabla\varphi
\end{align}
as $\epsilon=\epsilon_j\searrow 0$. In addition, a simple calculation with \eqref{lem332} shows that
\begin{align}\label{t116}
|-\epsilon\int_0^T\int_\Omega \frac{u_\epsilon^{k+1}}{1+u_\epsilon}\varphi|&\le \epsilon\|\varphi\|_{L^\infty(\Omega\times(0,T))}\int_0^T\int_\Omega u_\epsilon^k\nonumber\\
&\le\epsilon M_1(1+T)\|\varphi\|_{L^\infty(\Omega\times(0,T))}\rightarrow 0
\end{align}
as $\epsilon=\epsilon_j\searrow 0$. Consequently, we obtain from \eqref{t19} with \eqref{t110}--\eqref{t116} and the Fatou lemma that
\begin{align}\label{t120}
\int_0^T\int_\Omega \Big(\frac{\nabla u}{1+u}-&\frac{\chi u\nabla v}{2(1+u)v}\Big)^2\varphi\le\liminf_{\epsilon=\epsilon_j\searrow 0}\int_0^T\int_\Omega \Big(\frac{\nabla u_\epsilon}{1+u_\epsilon}-\frac{\chi u_\epsilon\nabla v_\epsilon}{2(1+u_\epsilon)v_\epsilon}\Big)^2\varphi \nonumber\\
 &=-\int_0^T\int_\Omega \ln(1+u)\varphi_t-\int_\Omega \ln(1+u_0)\varphi(\cdot,0)+\frac{\chi^2}{4}\int_0^T\int_\Omega\frac{u^2|\nabla v|^2}{(1+u)^2v^2}\varphi\nonumber\\
&~~+\int_0^T\int_\Omega\frac{\nabla u \cdot\nabla \varphi}{1+ u}-\chi\int_0^T\int_\Omega\frac{u}{(1+u)v} \nabla v\cdot\nabla\varphi\nonumber\\
&~~
-r\int_0^T\int_\Omega \frac{u}{1+u}\varphi+\mu\int_0^T\int_\Omega \frac{u^k}{1+u}\varphi
\end{align}
as $\epsilon=\epsilon_j\searrow 0$. This together with \eqref{t18} yields that $(u,v)$ is a weak logarithmic supersolution of \eqref{p} as well.\medskip

The proof is complete.
\qquad$\Box$\medskip

Next, we will prove that the global generalized solution to \eqref{p} is globally bounded. At first, we give a crucial estimate on $\int_\Omega u_\epsilon^pdx$ for all $\epsilon\in(0,1)$.
\begin{lemma}\label{lemma311}
Let $(u_\epsilon,v_\epsilon)$ be the global very weak solution of the problem \eqref{p} established in Theorem \ref{th2}. Then for $p>\frac{n(n+2)}{2(n+1)}$ we have
\begin{align}\label{lem311}
\frac{d}{dt}\int_\Omega u_\epsilon^pdx& \le -\int_\Omega u_\epsilon^pdx+M_6\Big(\int_\Omega u_\epsilon^pdx\Big)^\frac{q-qa}{p-qa}\nonumber\\
&~~+M_6\Big(\int_\Omega u_\epsilon^p dx\Big)^\frac{q}{p}+M_6\Big(\int_\Omega u_\epsilon^p dx\Big)^\frac{2q}{p(q-p)}+M_6\int_\Omega u_\epsilon dx,~t>0
\end{align}
for all $\epsilon\in(0,1)$ with some $M_6>0$.
\end{lemma}
{\bf Proof.} \
It follows from \eqref{pp}$_1$ and \eqref{lem401} for $1<p<q$ that
\begin{align}\label{3111}
\frac{1}{p}\frac{d}{dt}\int_\Omega u_\epsilon^pdx&=\int_\Omega u_\epsilon^{p-1}[\Delta u_\epsilon-\chi\nabla\cdot(\frac{u_\epsilon}{v_\epsilon}\nabla v_\epsilon)+r u_\epsilon-\mu  u_\epsilon^k-\epsilon u_\epsilon^{k+1}]dx\nonumber\\
&\le -\frac{1}{p}\int_\Omega u_\epsilon^pdx-(p-1)\int_\Omega u_\epsilon^{p-2}|\nabla u_\epsilon|^2dx+{\chi(p-1)}\int_\Omega \frac{u_\epsilon^{p-1}}{v_\epsilon^2}\nabla u_\epsilon\cdot\nabla v_\epsilon dx\nonumber\\
&~~+(r+1)\int_\Omega u_\epsilon^pdx-\mu \int_\Omega u_\epsilon^{p+k-1}dx\nonumber\\
&\le -\frac{1}{p}\int_\Omega u_\epsilon^pdx-\frac{p-1}{2}\int_\Omega u_\epsilon^{p-2}|\nabla u_\epsilon|^2dx+
\frac{\chi^2(p-1)}{2\delta_1^2}\int_\Omega u_\epsilon^p |\nabla v_\epsilon|^2dx\nonumber\\
&~~+(r+1)\int_\Omega u_\epsilon^pdx-\frac{\mu}{2}\int_\Omega u_\epsilon^{p+k-1}dx\nonumber\\
&\le -\frac{1}{p}\int_\Omega u_\epsilon^pdx-\frac{p-1}{2}\int_\Omega u_\epsilon^{p-2}|\nabla u_\epsilon|^2dx+\frac{\chi^2(p-1)}{2\delta_1^2}\int_\Omega u_\epsilon^qdx\nonumber\\
&~~+\frac{\chi^2(p-1)}{2\delta_1^2}\int_\Omega |\nabla v_\epsilon|^\frac{2q}{q-p}dx+C_{14}\int_\Omega u_\epsilon dx,~~t>0
\end{align}
by Young's inequality with $C_{14}=(r+2)^\frac{p+k-1}{k-1}(\frac{2}{\mu})^\frac{n}{k-1}$.
Invoking the Gagliardo-Nirenberg inequality, we get
\begin{align}\label{31111}
\|u_\epsilon\|_{L^q(\Omega)}&=\|u_\epsilon^\frac{p}{2}\|_{L^\frac{2q}{p}(\Omega)}^\frac{2}{p}\le C_{GN}\|u_\epsilon^\frac{p}{2}\|_{W^{1,2}(\Omega)}^\frac{2a}{p}\|u_\epsilon^\frac{p}{2}\|_{L^2(\Omega)}^{\frac{2(1-a)}{p}}\nonumber\\
&\le 2^\frac{pa}{2}C_{GN}\Big(\|u_\epsilon^\frac{p-2}{2}\nabla u\|_{L^2(\Omega)}^\frac{2a}{p}\|u_\epsilon^\frac{p}{2}\|_{L^2(\Omega)}^{\frac{2(1-a)}{p}}+\|u_\epsilon^\frac{p}{2}\|_{L^2(\Omega)}^{\frac{2}{p}}\Big).
\end{align}
If $1<p<q<\frac{n+2}{n}p$, we know  $a=\frac{n}{2}-\frac{pn}{2q}\in(0,1)$ and $\frac{2qa}{p}<1$. This fact together with \eqref{31111} yields
\begin{align}\label{31112}
\frac{\chi^2(p-1)}{2\delta_1^2}\int_\Omega u_\epsilon^qdx\le \frac{p-1}{2} \int_\Omega u_\epsilon^{p-2}|\nabla u_\epsilon|^2 dx+C_{15}\Big( \int_\Omega u_\epsilon^p dx\Big)^\frac{q(1-a)}{p-qa}+C_{16}\Big(\int_\Omega u_\epsilon^px\Big) ^\frac{q}{p}
\end{align}
by Young's inequality with $C_{15}=(p-1)^\frac{p}{qa}(2^{q+\frac{pqa}{2}}\frac{\chi^2}{\delta_1^2}C_{GN}^q)^\frac{p}{p-qa}$ and $C_{16}=2^{q+\frac{pqa}{2}}C_{GN}^q\frac{\chi^2}{\delta_1^2}(p-1)$.
Now, let $ p\in(\frac{n(n+2)}{2(n+1)},n]$ with $p<q<\frac{n+2}{n}p$. Then $\frac{2q}{q-p}<\frac{np}{n-p}$. By the classical imbedding Theorem with \eqref{443}, we obtain
\begin{align}\label{3112}
\frac{\chi^2(p-1)}{2\delta_1^2}\int_\Omega |\nabla v_\epsilon|^\frac{2q}{q-p}dx&=\frac{\chi^2(p-1)}{2\delta_1^2}\|\nabla v_\epsilon\|_{L^\frac{2q}{q-p}(\Omega)}^\frac{2q}{q-p}\nonumber\\
&\le C_{17}\frac{\chi^2(p-1)}{2\delta_1^2}\| \nabla v_\epsilon\|_{W^{1,p}(\Omega)}^\frac{2q}{q-p}\nonumber\\
&\le C_{18}\Big(\int_\Omega  u_\epsilon^pdx\Big)^\frac{2q}{p(q-p)},~t>0
\end{align}
with some  $C_{17},C_{18}>0$.
Combing \eqref{3111} with \eqref{31112} and \eqref{3112}, we have
\begin{align}\label{3113}
\frac{1}{p}\frac{d}{dt}\int_\Omega u_\epsilon^pdx& \le -\frac{1}{p}\int_\Omega u_\epsilon^pdx+C_{15}\Big(\int_\Omega u_\epsilon^pdx\Big)^\frac{q-qa}{p-qa}\nonumber\\
&~~+C_{16}\Big(\int_\Omega u_\epsilon^p dx\Big)^\frac{q}{p}+C_{18}\Big(\int_\Omega u_\epsilon^p dx\Big)^\frac{2q}{p(q-p)}+C_{14}\int_\Omega u_\epsilon dx,~t>0.
\end{align}
This completes the conclusion \eqref{lem311} with $M_6=p\max\{C_{14},C_{15},C_{16},C_{18}\}.$
\qquad$\Box$\medskip

Now, we establish the following uniform-in-time estimate on $\int_\Omega u_\epsilon^pdx$ for all $\epsilon\in(0,1)$, with the initial data $u_0$ and $\frac{r}{\mu}$ suitably small.
\begin{lemma}\label{lemma312}
Let $(u_\epsilon,v_\epsilon)$ be the global very weak solution of the problem \eqref{p} established in Theorem \ref{th2}. Then for $p\in(\frac{n(n+2)}{2(n+1)},n]$ there exist $\eta,\lambda>0$ such that
\begin{align}\label{lem312}
\int_\Omega u_\epsilon^pdx\le M_7,~~t>0,
\end{align}
provided $\frac{r}{\mu}<\eta$ and $\int_\Omega u_0^pdx<\lambda$, for all $\epsilon\in(0,1)$ with    $M_7>0$.
\end{lemma}
{\bf Proof.}\
Let $F_\epsilon(t):=\int_\Omega u_\epsilon(x,t)^pdx,~t>0$. Then we have from \eqref{3113} and \eqref{lem331} that
\begin{equation}\label{3115}
\begin{cases}
F_\epsilon'(t)\le -F_\epsilon(t)+M_6F_\epsilon(t)^\frac{q-qa}{p-qa}+M_{6}F_\epsilon(t)^\frac{q}{p}+M_6F_\epsilon(t)^\frac{2q}{p(q-p)}+M_{6}m^*,~~t>0,\\
F_\epsilon(0)=\int_\Omega u_0^pdx.
\end{cases}
\end{equation}
Since $p\in(\frac{n(n+2)}{2(n+1)},n]$ and $p<q<\frac{n+2}{n}p$, we know  $\frac{q-qa}{p-qa},\frac{q}{p},\frac{2q}{p(q-a)}>1$. Denote $$h(s,m^*):=-s+M_6s^\frac{q-qa}{p-qa}+M_{6}s^\frac{q}{p}+M_{6}s^\frac{2q}{p(q-p)}+M_{6}m^*,~~s>0.$$
Then there exists $m_0^*>0$ such that $h(s,m_0^*)$ has the unique positive root $s_0$. Furthermore,
$M(t)\equiv s_0$ verifies the ODE problem
\begin{equation}\label{3116}
\begin{cases}
M'(t)=h(M(t),m_0^*),~~t>0,\\
M(0)=s_0.
\end{cases}
\end{equation}
If $m^*<m_0^*$,  it follows by a continuous dependence argument that the function
$h(s,m^*)$, with $h(s,m^*)<h(s,m_0^*)$,  has exactly two positive roots $0<s_1<s_0<s_2$. Now let
\begin{align*}
\eta:=\Big(\frac{m_0^*}{|\Omega|}\Big)^{k-1}~~{\rm and~~}\lambda:=\min\Big\{s_0,\frac{m_0^{*p}}{|\Omega|^{p-1}}\Big\}
\end{align*}
with $\frac{r}{\mu}<\eta$ and $\int_\Omega u_0^pdx<\lambda$. Then
$$\int_\Omega u_0dx<|\Omega|^\frac{p-1}{p}\Big(\int_\Omega u_0^pdx\Big)^\frac{1}{p}<m_0^*~~{\rm and}~~\int_\Omega u_\epsilon dx\le\max\Big\{\int_\Omega u_0dx,(\frac{r}{\mu})^\frac{1}{k-1}|\Omega|\Big\}<m_0^*$$
for all $\epsilon\in(0,1)$.
 Consequently, we obtain from these estimates with problems \eqref{3115} and \eqref{3116} that
$$F_\epsilon(t)=\int_\Omega u_\epsilon^pdx\le s_1,~~t>0$$
for all $\epsilon\in(0,1)$ by an ODE  comparison principle.
 The proof is complete. \qquad$\Box$\medskip\\
{\bf Proof of Theorem \ref{th3}}\
Based on the estimate $\int_\Omega u_\epsilon ^pdx\le s_1$ for $p>\frac{n(n+2)}{2(n+1)}$ in Lemma \ref{lemma311} and uniformly in time lower-bound estimate of $v_\epsilon$, we obtain the global boundedness of solutions to the regularization problem \eqref{pp} via a similar argument as that in \cite[Lemma 2.3]{M3}, i.e., $\|u_\epsilon\|_{L^\infty(\Omega)}\le \tilde{C}$ with some $\tilde{C}>0$ for all $t>0$ and $\epsilon\in(0,1)$. Consequently, we conclude that the generalized solution $(u,v)$ is globally bounded as well by taking  $\epsilon=\epsilon_j\searrow 0$. \qquad$\Box$\medskip

{\noindent \bf Acknowledgements}

This work was supposed by the Doctoral Scientific Research
Foundation of Liaoning Normal University (Grant No. 203070091907).



{\small \begin{thebibliography}{}

\bibitem{EL} E. F. Keller, L. A. Segel, Initiation of slime mold aggregation viewed as an instability, J. Theoret. Biol. 26 (1970) 399--415.\\[-17pt]
\bibitem{KMT} K. Fujie, M. Winkler, T. Yokota, Boundedness of solutions to parabolic-elliptic Keller-Segel systems with signal-dependent sensitivity, Math. Methods. Appl. Sci. 38 (2015) 1212--1224.\\[-17pt]
\bibitem{KOA} K. Osaki, A. Yagi, Finite dimensional attractors for one-dimensional Keller-Segel equations, Funkcial Ekvac. 44 (2001) 441--469.\\[-17pt]
\bibitem{FS} K. Fujie, T. Senba, Global existence and boundedness in a parabolic-elliptic Keller-Segel
  system with general sensitivity, Discrete Contin. Dyn. Syst. Ser. B 21 (2016) 81--102.\\[-17pt].
\bibitem{B} T. Black, Global generalized solutions to a parabolic-elliptic Keller-Segel system with singular sensitivity, Discrete Contin. Dyn. Syst. Ser. S 13 (2020) 119--137.\\[-17pt]
\bibitem{K} K. Fujie, Boundedness in a fully parabolic chemotaxis system with singular sensitivity, J. Math. Anal. Appl. 424 (2015) 675--684.\\[-17pt]
\bibitem{FST} K. Fujie, T. Senba, Global existence and boundedness of radial solution to a two dimensional fully parabolic chemotaxis system with general sensitivity, Nonlinearity 29 (2016) 2417--2450.\\[-17pt]
\bibitem{LW} J. Lankeit, M. Winkler, A generalized solution concept for the Keller-Segel system with logarithmic sensitivity: global solvability for large nonradial data, Nonlinear Differ. Equ. Appl. (2017) 24--49.\\[-17pt]
\bibitem{TT} T. Nagai, T. Senba, Behavior of radially symmetric solutions of a system related to chemotaxis, Nonlinear Anal. 30 (1997) 3837--3842.\\[-17pt]
\bibitem{JM} J. I. Tello, M. Winkler, A chemotaxis system with logistic source, Comm. Partial
Differential  Equations 32 (2007) 849--877. \\[-17pt]
\bibitem{M3} M. Winkler, Chemotaxis with logistic source: Very weak global solutions
and their boundedness properties, J. Math. Anal. Appl. 348 (2008) 708--729.\\[-17pt]
\bibitem{M2} M. Winkler, Blow-up in a higher-dimensional chemotaxis system despite logistic
growth restriction, J. Math. Anal. Appl. 384 (2011) 261--272.\\[-17pt]
\bibitem{OO}    K. Osaki, T.Tsujikawa, A. Yagi, M. Mimura, Exponential attractor for a chemotaxis-growth
system of equations, Nonlinear Anal. 51 (2002) 119--144. \\[-17pt]
\bibitem{W1} M. Winkler,
  Boundedness in the higher-dimensional parabolic-parabolic chemotaxis system with logistic source,
 Comm. Partial  Differential  Equqtions  35 (2010) 1516--1537. \\[-17pt]
\bibitem{V} G. Viglialoro, Very weak global solutions to a parabolic-parabolic
chemotaxis-system with logistic source, J. Math. Anal. Appl. 439 (2016) 197--212.\\[-17pt]
\bibitem{V1}  G. Viglialoro, Boundedness properties of very weak solutions to a fully parabolic
chemotaxis-system with logistic source, Nonlinear Anal. RWA 34 (2017) 520--535.\\[-17pt]
\bibitem{L1} J. Lankeit, Eventual smoothness and asymptotics in a three-dimensional chemotaxis system with logistic source, Journal of Differential Equations 258 (2015) 1158--1191.\\[-17pt]
\bibitem{TM} Y. S. Tao, M. Winkler, Persistence of mass in a chemotaxis system with logistic source,
J. Differential Equations 259 (2015) 6142--6161. \\[-17pt]
\bibitem{M4} M. Winkler, How strong singularities can be regularized by logistic degradation in the Keller-Segel system? Ann. Mat. Pura Appl. 198 (2019) 1615--1637.\\[-17pt]
\bibitem{FF} K. Fujie, M. Winkler, T. Yokota, Blow-up prevention by logistic sources in a parabolic-elliptic Keller-Segel system with singular sensitivity, Nonliear Anal. 109 (2014) 56--71. \\[-17pt]
\bibitem{ZZ} X. D. Zhao, S. N. Zheng, Global boundedness to a chemotaxis system
with singular sensitivity and logistic source, Z. Angew. Math. Phys. 68:2 (2017) 13 pp.
\bibitem{ZZZ} X. D. Zhao, S. N. Zheng, Global existence and boundedness of solutions to
a chemotaxis system with singular sensitivity and
logistic-type source, J. Differential Equations  267 (2019) 826--865.\\[-17pt]
\bibitem{LL1} E. Lankeit, J. Lankeit, Classical solutions to a logistic chemotaxis model with singular sensitivity and signal absorption,
Nonlinear Anal. RWA 46 (2019) 421--445.\\[-17pt]
\bibitem{LL} E. Lankeit, J. Lankeit, On the global generalized solvability of a chemotaxis model with signal absorption and logistic
growth terms, Nonlinearity 32 (2019) 1569--1596.\\[-17pt]
\bibitem{Z1}X. D. Zhao, S. N. Zheng, Global existence and asymptotic behavior to a chemotaxis¨Cconsumption system with singular sensitivity and logistic source, Nonlinear Anal. RWA 42 (2018) 120--139. \\[-17pt]
\bibitem{M5} M. Winkler, The two-dimensional Keller-Segel system with singular sensitivity and signal absorption: Global large-data solutions and their relaxation properties, Math. Models Methods Appl. Sci. 26(5) (2016) 987--1024.\\[-17pt]
\bibitem{M6} M. Winkler, The two-dimensional Keller-Segel system with singular sensitivity and signal absorption: Eventual smoothness and equilibration of small-mass solutions, 2016, preprint.
\bibitem{M7} M. Winkler, Renormalized radial large-data solutions to the higher-dimensional Keller-Segel system with singular sensitivity and signal absorption, J. Differential Equations 264 (2018) 2310--2350. \\[-17pt]
\bibitem{SSW} C. H. Stinner, C. H. Surulescu, M. Winkler, Global weak solutions in a PDE-ODE
system modeling multiscale cancer cell invasion, SIAM J. Math. Anal. 46 (2014) 1969--2007. \\[-17pt]
\bibitem{GST} E. Galakhov, O. Salieva, J. I. Tello,
On a parabolic-elliptic system with chemotaxis and logistic type growth, J. Differential Equations 261 (2016) 4631--4647. \\[-17pt]
\bibitem{BS}H. Br\'{e}zis, W. A. Strauss, Semi-linear second-order elliptic equations in $L^1$, J. Math. Soc. Japan 25 (1973) 565--590.\\[-17pt]
 \bibitem{M1} M. Winkler, Aggregation vs. global diffusive behavior in the higher-dimensional Keller-Segel model,
 J. Differential Equations 248 (2010) 2889--2905.\\[-17pt]


\end{thebibliography}}
\end{document}